# A mathematical model for planning oil products distribution via pipeline

Amir Baghban, Department of Mathematics and Computer Sciences, Amirkabir University of Technology, Tehran, Iran

Majid Yousefikhoshbakht, Department of Mathematics, Faculty of Sciences, Bu Ali Sina University, Hamedan, Iran

**Abstract** Compared to other transportation modes, e.g., road, railroad and vessel, multiproduct pipelines are the safest and most economical way of conveying petroleum products over long distances day and night. During the last two decades, the operational multiproduct pipeline scheduling has gained increasing attention, where most of the contributions are based on continuous-time representation. In this paper, we present a new discrete-time mixed-integer linear programming (MILP) model for the short-term scheduling of multiproduct pipelines featuring multiple refineries and distribution centers. The model optimally determines the sequence/ planning of product injections at input nodes, the sequence of product deliveries to output nodes, and traces the size and the product's location along the pipeline at any time. The proposed model also allows to execute simultaneous injections at the input nodes and to tackle simultaneous injections and deliveries at an intermediate node that can act as both the input and output node (i.e., dual-purpose nodes). Solutions to two benchmark example problems illustrate that the proposed model presents significant reductions in the pipeline operational cost.

## 1. INTRODUCTION

In comparison with different transportation modes, e.g., train, truck and vessel, the pipelines are more safe and economical and available every time. Suffice it to say that the oil pipelines account for approximately 70% of all refined oil products shipment in the US [4]. A plethora of works can be found in the supply chain optimization problems that focus on reducing carbon emissions [1, 19, 20, 28,45]. In terms of downstream oil supply chain, the oil pipelines generate roughly no environmental pollution and can transport diverse oil products (e.g. heating oil, gasoline, kerosene, fuel of jets, liquid petroleum gases, diesel etc.).

As an operational rule, the volume of oil products inside the pipeline must be equal to the total volume of pipeline at any time, and consequently, when some volume of oil products is injected into the pipeline, the same amount should be removed at the other end. Since products are moving forward along the pipeline without physical seperation, there is always some contaminated volume of adjacent products. The volume of the mixing volume (known for any pipeline) is called *interface* in the oil pipeline industry. The volume usually needs to be reprocessed at a refinery. The degree of the losses due to the mixing volume depends on the type of two products that move successively. If the interface of the two products is large, the decision makers do not allow them to be injected successively. For instance, gas oil should never be pumped into the pipeline after/before gasoline. With a proper sequence of oil products in a pipeline, the oil pipeline company can significantly reduce the interface costs.

The operational scheduling of pipeline systems is a complex optimization task such that many constraints should be included. Constraints such as (1) monitoring product inventories at input and output nodes, (2) tracking the size and the location of products during the time horizon at time points, and (3) respecting forbidden product sequences are only a short list of limitations [18]. The main challenge is to provide the right product to the right customer until the end of planning horizon at the lower operational cost (e.g. pumping, interface, and inventory costs).

The pipeline scheduling problem (PSP) has been studied using different types of optimization approaches including heuristic techniques ([2,3,16, 22,30-31,39]), decomposition frameworks ([27,29]) and MILP formulations [5,14,40,43,47-49]. The latter can be divided into two categories: discrete and continuous-time. In discrete approaches, both of the planning horizon and the pipeline length are respectively divided into time slots of fixed duration and into lots of uniform sizes[17, 26, 42-43]. These limitations are relaxed in the continuous methods and the time slot length is determined by optimization [5, 10, 11, 15, 31].

Previous works on the PSP have considered pipeline systems with a variety of configurations involving straights, tree, and mesh structures [24, 27, 32, 40,

46,48]. Cafaro and Cerda [4, 5] were the first to present a continuous-time model for scheduling aggregately a pipeline comprised of a single refinery and multiple output nodes. Later they extended their approach in [4] to study the scheduling of multi-refinery [6, 7], tree-like [8], and mesh structured pipeline systems [9]. Based on Cafaro's model [4], Relvas et al. [43] studied the scheduling of a multiproduct pipeline featuring a single refinery and just one terminal to consider inventory management at depots. Mostafaei and coworkers [13, 14, 21, 33-38,46] used continuous-time MILP models to study the detailed scheduling and lot-sizing problems of the pipeline systems with straight and tree-like structures. Liao et al. [23-24] proposed MILP models that allow to inject multiple products during a time window, thus resulting in good schedules in lower computational times.

Rejowski and Pinto [42] proposed the first discrete-time model for the PSP. The approach considers the problem of pipeline scheduling for a pipeline network with just one refinery and multiple terminals. The model meets all the operational constraints such as mass balances, product demands, and delivery and constraints on inventory. Later, they improved their previous work in [42] by applying integer cuts and imposing extra practical constraints to overcome the problem complexity [43]. Note that the models by Rejowski and Pinto [42-43] assume that the product injection rate is a constant number.

To the best of our knowledge, all contributions on the scheduling of pipeline systems with multiple inputs and output nodes, also the subject of this paper, have been developed based on continuous-time MILP models [6-7,12,14,24,36,44]. Note that, as in continuous-time models the number of time slots is difficult to predict a priori, the user needs an iterative procedure to find the optimal solutions. Therefore, a particular example problem needs to be solved at least two times [6-7,24].

Discrete-time models, on the other hand, can solve a particular problem in a single iteration [34, 42]. The main contribution of this article is then to develop a discrete-time MILP model for the pipeline scheduling of pipelines with multiple inputs and output nodes. The proposed model allows tackling the pumping and delivery operations at the dual purpose nodes that can be performed simultaneously. It is also worth mentioning that the existing discrete-time scheduling models [17, 26, 34, 38, 41-43] can only handle the operation of the straight pipeline systems with a single input node and multiple distribution centers, and some of them assume a constant pumping rate at the refinery. The Table 2 compares some articles (from the earliest ones to the recently published ones) with the current article in terms of the topology of the pipeline and the efficiency of the models developed in them.

In this paper, we propose a discrete-time-based MILP model that can handle variable pumping rates at the input nodes. The goal is to fulfill the demands of products at the output nodes at the lowest pumping, interface, and backorder costs while optimally determining: (1) the sequence and timing of product injections at the input nodes, (2) the size of products inside the pipeline, and (3) the sequence and timing of product deliveries to depots. The rest of the paper is structured as follows: Section 2 illustrates an example explaining the efficiency of the proposed model and Section 3 explains the pipeline structure and main assumptions supposed in the model. The next section develops a discrete-time MILP model for a straight pipeline that connects multiple input nodes to multiple output nodes. In section 5, we benchmark the proposed optimization approach with a continuous-time model from the literature. Finally, Section 6 gives the conclusions and some future directions.

| NOMENCLATURE | | | |
|---|---|---|---|
| Sets | | $\max_p$ | Upper bound for injection of product $p$ during the time horizon |
| $T$ | Set of time slots indicated by $t$, $t = 1, \ldots, \|T\|$ | $S_{\min}, S_{\max}$ | Minimum and maximum amount of materials that can move from line $l$ to $l-1$ |
| $I$ | Product batches indexed by $i$, $i = 1, \ldots, \|I\|$ | $W_{i,l}^{initial}$ | A parameter indicating the volume of old batch $i$ in line $l$ at the start time of the planning horizon |
| $I_l^{old}$ | Old batches that fill line $l$ at the start of planning horizon | $Cp_{p,l}$ | Pumping cost per unit of product $p$ by $l$ |
| $I^{new}$ | New batches | $cb$ | Backorder cost |
| $L$ | Set of pipeline segments placed between two successive refineries indexed by l | $Interface_{p,p'}$ | Interface cost related to product sequence $(p,p')$ |
| $D_l$ | Set of distribution centers after line $l$ and before line $l+1$ indexed by $d$, $d = 1, \ldots, |D_l|$ | Continuous variables | |
| $I_l^{input}$ | Set of batches that refinery $l$ can inject | $F_{i,l,t}$ | The upper volumetric coordinate of batch $i$ in line $l$ with respect to refinery $l$ during the time slot $t$ |
| $I_l$ | Batches that can be pushed forward in line $l$ | $W_{i,l,t}$ | The volume of batch $i$ contained in line $l$ at time slot $t$ |
| $P$ | Set of oil products that are supposed to be pumped being indexed by $p$, $p = 1, \ldots, |P|$ | $R_{i,l,t}$ | The volume of batch $i$ that input node $l$ injects during time slot $t$ |
| Parameters | | $RP_{i,p,l,t}$ | Injected volume of product $p$ into the batch $i$ by the refinery $l$ during time slot $t$ |
| $h_{max}$ | The length of the planning horizon in hours (h) | $D_{i,d,l,t}$ | Size of batch $i$ received by depot $d$ located on line $l$ during time slot $t$ |
| $FP_{p,p'}$ | A Boolean matrix indicating the possible sequence of successive batches between products $(p,p')$; $FP_{p,p'} = 1$ denotes that product $p$ can be adjacent to $p'$. | $DP_{i,p,d,l,t}$ | The volume of product $p$ delivered to depot $d$ from batch $i$ in line $l$ during time slot $t$ |
| $f_{i,l}^{initial}$ | A parameter specifying the upper coordinate of old batch $i$ in line $l$ | $Border_{p,d,l}$ | Backorder volume of product $p$ corresponding to depot $d$ in line $l$ |
| $vl_l$ | The volume of line $l$ in cubic meter | $S_{i,l,t}$ | The volume of the transferred batch $i$ from line $l-1$ to line $l$ during time slot $t$ |
| $r_l^{min}, r_l^{max}$ | Lower and upper bounds for injection operation at the refinery $l$ during a time slot | $Ic_i$ | Reprocessing cost spent between batches $i$ and $i-1$ |
| $vr_l^{min}, vr_l^{max}$ | Upper and lower bounds for injection rate at node $l$ | Binary variables | |
| $\Delta t$ | Length of slot times | $y_{i,p}$ | The assignment variable for product $p$ to batch $i$ |
| $invp_{p,l}$ | The amount of product $p$ existing at the tanks of refinery $l$ at $t = 0$ | $w_{i,l,t}$ | 1 if the input node $l$ injects the batch $i$ during time interval $t$ |
| $\sigma_d$ | The distance of output node $d$ from $l$, located on line $l$ | $x_{i,d,l,t}$ | 1 if the distribution center $d$, on line $l$, receives batch $i$ in time slot $t$ |
| $d_{min}, d_{max}$ | Minimum and maximum amounts of a batch that can be received by an output node | $u_{i,l,t}$ | 1 if batch $i$ is transferred from line $l-1$ to line $l$ at time slot $t$ |
| $demand_{p,d,l}$ | The volume of product $p$ ordered by output $d_l$ | | |

**Table 1.** comparisons between the current article and existing ones

| Reference number | One refinery and output node | One refinery and multiple outputs | multiple inputs and multiple outputs | Discrete time | Continuous time | Dual purpose station | Simultaneous injection | Variable pumping rate |
|---|---|---|---|---|---|---|---|---|
| [38] (2003) | | * | | * | | | | |
| [41] (2003) | | * | | * | | | | |
| [5] (2004) | | * | | | * | | | |
| [26] (2004) | * | | | * | | | | * |
| [6] (2010) | | | * | | * | * | * | * |
| [10] (2015) | | | * | | * | * | * | * |
| [36] (2016) | | | * | | * | * | * | * |
| [17] (2017) | | * | | * | | | | * |
| [30] (2018) | | * | | | * | | | * |
| [37] (2021) | * | * | * | | * | * | * | * |
| [34] (2021) | * | * | | * | | | | * |
| This work (2021) | * | * | * | * | | * | * | * |

## 2. MOTIVATING EXAMPLE

The scheduling problem of multiproduct pipelines is an NP-hard problem and the mathematical model of them consists of many continuous and binary variables to be specified. These models are able to allow the refineries to perform simultaneous injections at variable pumping rates during the planning horizon that contributes to minimizing the makespan and reducing the operational costs (due to reduction in some pipeline segments being idle). The majority of operational costs are due to the idle state of some parts of the pipelines. However, developed discrete-time models are only capable of handling pipelines with single input nodes. So, the main novelty of this article is scheduling operational activities of pipelines with multiple input and output nodes using a discrete-time mathematical model, where it is able to schedule simultaneous injections at variable pumping rates.

In this section, we will illustrate, by an example, how cost-effective it would be to incorporate constraints handling simultaneous injections with variable pumping rate in the mathematical model. The pipeline of this example involves two refineries (L1 and L2) with three output nodes (D1, D2, D3). The initial linefill for the left half of pipeline is comprised of 40 units of product P2 in batch I2 and the second half part of the pipeline contains batch I1 conveying 40 units of product P1. This example tries to plan the delivery of 22 units of product P2 to depot D1 in line L1 and 5 units of product P1 to each of depots D1 and D2 in line L2 in just 20 hours.

As depicted in Figure 1, if simultaneous injections and variable pumping rates are not allowed, there would be a delay of 20 hours in fulfilling the demands in output nodes and it can meet product demand after 4 time slots (it is worth noticing that each time slot duration is 10 hours). Moreover, due to the constant pumping rate, refinery L1 has to pump 8 extra units of P2 into the line, imposing extra pumping cost. However, allowing simultaneous injections with variable pumping rates can reduce the time spent to fulfill the demands by half (see Figure 2) so incredible reduction in the pumping costs can be seen.

## 3. PROBLEM STATEMENT

This article aims at scheduling straight pipelines with |L| input nodes and |D| output nodes where some of the facilities are dual purpose, i.e., capable of pumping products into the pipeline and receiving products from the pipeline. Regarding the lines and segments, parts of pipeline lying between two consecutive nodes are called lines numbered with the same indices as that of the left-hand refinery. Each line includes segments ending with outputs. For instance, Figure 3 can be an example of a pipeline of this type. It has two lines, meaning that the pipeline has two refineries and also it has three distribution centers. With regards to the latter, two of them can only receive products (D1 and D1') and the other can inject and receive products i.e., it is dual purpose output node (The node indicated by both D2 and L2).

Mathematical models to schedule the transportation of the oil products are either batch-centric such as the proposed model in this article or product centric. In batch centric models, the products are assigned to the batches with the smallest positioned at the right extreme of the pipeline (I1). In their sequencing, batch $I_{n+1}$ always comes after batch $I_n$. The batches contained in the initial linefill are called old batches where their contents and volumes are specified by the decision maker, while new batches are assigned by the model.

The following predefined conditions will be taken into account in the pipeline formulation:

(A1) Pipelines are always full of product and products will flow from left to the right unidirectionally (the products can not move reversely).
(A2) The amount of time spent on switching from one tank to another is neglected.
(A3) The volume of refined product mixture at the interface points is neglected.
(A4) There are no other sources than refineries adding products to tanks during the planning horizon.
(A5) Products are assumed to be incompressible and so the injected volume during each time slot will be discharged at the same time.

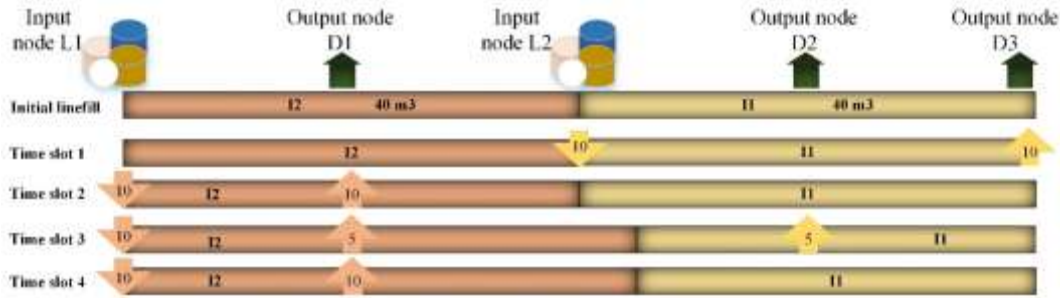

**FIGURE 1.** Pipeline scheduling without simultaneous injections with constant pumping rate

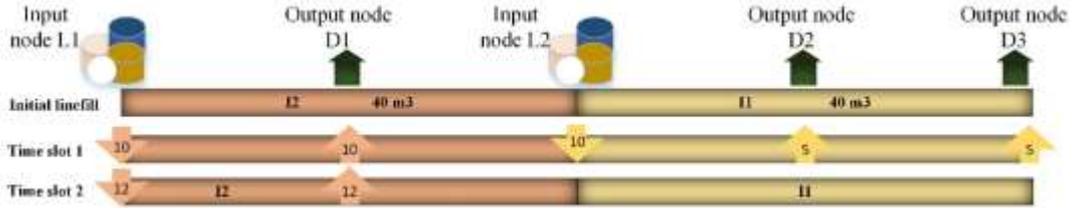

**FIGURE 2.** Pipeline scheduling by a mathematical model employing simultaneous injections and varying pumping rate

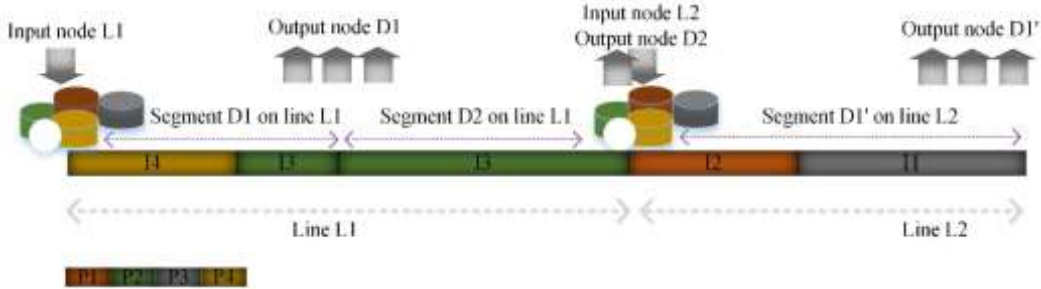

**FIGURE 3.** Illustration of a pipeline structure

### 4. MATHEMATICAL MODEL
This section is dedicated to propose a discrete-time MILP model to schedule a straight multi product pipeline connecting multiple input nodes to multiple output nodes. The model is inspired by continuous-time model developed in [37].

**4.1. Allocating Batches to Products and Forbidden Product Sequence Constraints** According to Equation (1), each batch can only convey one product and the assignment is done by the binary variables $y_{i,p}$.

Note that the assignment of old batches $i \in I_l^{old}$ is known as parameter (i.e., $y_{i,p} = 1$, $\forall i \in I_l^{old}$ for some $p$). To avoid generating large interfaces inside the pipeline, the forbidden products $p$ and $p'$ should not flow successively inside the pipeline, as imposed by Equation (2).

**4.2. Controlling the Location of Batches** Let the continuous variable $F_{i,l,t}$ indicates the upper location of batch $i$ in line $l$ at time point $t$. Since the batch $i + 1$ always moves after batch $i$, $F_{i,l,t}$ will be equal to the amount of batch $i$ in the same line ($W_{i,l,t}$) plus the upper coordinate of the batch $i + 1$ in line $l$ (Equation 3). The upper position of $i \in I_l^{old}$ at the start time of scheduling horizon ($t = 0$) is given beforehand and equals to $f_{i,l}^{initial}$ (Equation 4). The product batches can move only in one direction in the pipeline (from input nodes to output nodes), as imposed by Equation (5). $F_{i,l,t}$ should not exceed the total volume of line $l$, as imposed by Equation (6).

**4.3. Injection at Input Nodes** Let the binary variable $w_{i,l,t}$ takes one if the input node $l$ enlarges batch $i$ at the time slot $t$. Equation (7) guarantees injection of at most single batch by every refinery during each time interval. At each time point, at least one of the input nodes must be active, as imposed by Equation (8). Moreover, a new batch $i \in I^{new}$ assigned to a product must be injected into the pipeline during some time slots (Equation (9)).

$$\sum_{p \in P} y_{i,p} = 1, \quad \forall i \in I \quad (1) \tag{1}$$

$$y_{i-1,p} + y_{i,p'} \leq 1 \quad \forall i \in I, p, p' \in FP_{p,p'} \tag{2}$$

$$F_{i,l,t} = W_{i,l,t} + F_{i+1,l,t}, \quad \forall i \in I_l, l, t \tag{3}$$

$$F_{i,l,t} = f_{i,l}^{initial}, \quad \forall i \in I_l, l, t = 0 \tag{4}$$

$$F_{i,l,t-1} \leq F_{i,l,t}, \quad \forall i \in I_l, l, t \tag{5}$$

$$F_{i,l,t} \leq vl_l, \quad \forall i \in I_l, l, t \tag{6}$$

$$\sum_{i \in I_l^{Input}} w_{i,l,t} \leq 1, \quad \forall l, t \neq |T| \tag{7}$$

$$\sum_l \sum_{i \in I_l^{Input}} w_{i,l,t} \geq 1, \forall t \, t \neq |T| \tag{8}$$

$$\sum_p y_{i,p} \leq \sum_t \sum_l w_{i,l,t}, \quad \forall i \in I^{new} \tag{9}$$

$$F_{i,l,t} - W_{i,l,t} \leq vl_l(1 - w_{i,l,t}), \quad \forall i \in I_l^{Input}, l, t \neq |T| \tag{10}$$

$$F_{i,l-1,t} \geq vl_{l-1} w_{i,l,t}, \quad \forall i \in I_l^{Input}, l \neq 1, t \neq |T| \tag{11}$$

$$w_{i,l,t} r_l^{min} \leq R_{i,l,t} \leq r_l^{max} w_{i,l,t}, \quad \forall i \in I_l^{Input}, l, t \neq |T| \tag{12}$$

$$\Delta t \times vr_l^{min} \sum_{i \in I_l^{Input}} w_{i,l,t} \leq \sum_{i \in I_l^{Input}} R_{i,l,t} \leq \Delta t \times vr_l^{max}, \quad \forall l, t \neq |T| \tag{13}$$

$$\sum_p RP_{i,p,l,t} = R_{i,l,t}, \quad \forall i \in I_l^{Input}, l, t \neq |T| \tag{14}$$

At time point $t$, input node $l$ can enlarge the batch $i \in I_l^{Input}$ if (i) the lower coordinate of the batch is located at the beginning origin of line $l$, (ii) and the upper coordinate of the batch in line $l - 1$ is equal to the total volume of line $l - 1$. This means that $w_{i,l,t} = 1$ should imply $F_{i,l,t} - W_{i,l,t} = 0$ and $F_{i,l-1,t} = vl_{l-1}$ (Equation (10) and Equation (11)). The volume injected to the line $l$ during the time interval $t$ by an active input node ($R_{i,l,t}$) should belong to a feasible region $[r_l^{min}, r_l^{max}]$ (Equation (12)).

$$\sum_{t \neq |T|} \sum_l RP_{i,p,l,t} \leq \max_p y_{i,p}, \quad \forall i, p \tag{15}$$

$$\sum_{t \neq |T|} \sum_{i \in I_l^{Input}} RP_{i,p,l,t} \leq invp_{p,l}, \quad \forall p, l \tag{16}$$

The pump rate at an active input node during any time interval (the injected volume divided by time slot length) should be maintained within an acceptable rate $[vr_l^{min}, vr_l^{max}]$ (Equation (13)).

Since only one product is assigned to each batch, the volume of product $p$ injected into the batch $i$ by input node $l$ during the interval $t$ ($R_{i,p,l,t}$), amounts to $R_{i,l,t}$, as illustrated by Equation (14). If batch $i$ does not contain product $p$, $R_{i,p,l,t}$ will be zero (Equation (15)). Moreover, the total volume of product $p$ injected from an input node $l$ during the scheduling horizon ($\sum_{t \neq |T|} \sum_{i \in I_l^{Input}} RP_{i,p,l,t}$) should not surpass the inventory of product $p$, i.e, $invp_{p,l}$ (Equation (16)).

**4.4. Delivery to Depots (output nodes)** Unlike the input nodes, output nodes can receive multiple batches during a time interval $t$. The binary variable $x_{i,d,l,t}$ indicates the discharging of batch $i$ by depot $d$ positioned on line $l$ at time point $t$. A batch $i$ can be received by a depot $d_l$, if the lower and upper coordinates of the batch satisfy Equation (17) and Equation (18), where $\sigma_{d,l}$ is the volumetric distance of depot $d$ from refinery $l$ ($x_{i,d,l,t} = 1 \Rightarrow F_{i,l,t-1} - W_{i,l,t-1} \leq \sigma_{d,l} \leq F_{i,l,t}$). Equation (19) shows that each depot can only receive feasible amount belonging to range $[d_{min}, d_{max}]$ and it is restricted by Equation (20).

Since each batch can transport only one of the products, the volume of product $p$ in batch $i$ transferred to output node $d_l$ at time interval t, which is equal to $DP_{i,p,d,l,t}$, will amount to $D_{i,d,l,t}$, as shown by Equation (22).

$$F_{i,l,t} - W_{i,l,t} \leq \sigma_{d,l} + (vl_l - \sigma_{d,l})(1 - x_{i,d,l,t}), \forall\, i \in I_{d,l}, d \in D_l, l, t\, (t \neq |T|) \quad (17)$$

$$F_{i,l,t} \geq x_{i,d,l,t}\sigma_{d,l}, \quad \forall\, i \in I_{d,l}, d \in D_l, l, t\, (t \neq |T|) \quad (18)$$

$$d_{\min}\, x_{i,d,l,t} \leq D_{i,d,l,t} \leq d_{\max}\, x_{i,d,l,t}, \quad \forall\, i \in I_{d,l},\, d \in D_l, l, t\, (t \neq |T|) \quad (19)$$

$$\sum_{d' \in D_l}^{d} D_{i,d,l,t} \leq \sigma_{d,l} - (F_{i,l,t} - W_{i,l,t}) + R_{i,l,t} + S_{i,l,t} + (vl_l - \sigma_{d,l})(1 - x_{i,d,l,t}), \forall\, i \in I_{d,l}, d \in D_l, l, t\, (t \neq |T|) \quad (20)$$

$$\sum_{\substack{t \\ t \neq |T|}} \sum_l \sum_{d \in D_l} DP_{i,p,d,l,t} \leq \max_p y_{i,p}, \quad \forall\, i, p \quad (21)$$

$$\sum_p DP_{i,p,d,l,t} = D_{i,d,l,t}, \forall\, i \in I_{d,l}, d \in D_l, l, t\, (t \neq |T|) \quad (22)$$

$$\sum_{\substack{t \\ t \neq |T|}} \sum_{i \in I_{d,l}} DP_{i,p,d,l,t} + Border_{p,d,l} \geq Demand_{p,d,l}, \quad \forall\, p, d \in D_l, l \quad (23)$$

If the batch does not contain product p, variable $DP_{i,p,d,l,t}$ will take zero (Equation (21)). To satisfy the product demands of the depots, the overal amount of product p delivered to this depot ($\sum_{\substack{t \\ t \neq |T|}} \sum_l \sum_{d \in D_l} DP_{i,p,d,l,t}$) should be equal to the model parameter $Demand_{p,d,l}$, otherwise, the backorder material will take a positive value in Equation (23) and results in backorder costs. Note that the continuous variable $Border_{p,d,l}$, standing for the backorder of product $p$ at the depot $d$ during the scheduling horizon will be minimized in the objective function.

**4.5. Transferring Material between adjacent lines**
Similarly to the depot operations, arbitrary number of batches can enter from one line to another during a time interval. Let the binary variable $u_{i,l,t}$ indicates moving of batch $i$ from segment $l - 1$ to segment $l$ during time slot $t$. To prevent products from being mixed, we do not allow batches to pass from line $l - 1$ to line $l$ during time interval $t$, if the input node $l$ is pumping some material into line $l$ over the same time slot by applying Equation (27). If $u_{i,l,t} = 1$, (1) then $F_{i,l,t}$ is at the end of line $l$ (Equation (25)), and (2) $F_{i,l,t} - W_{i,l,t}$ is equal to zero (Equation (24)). The transferred amount of product, from one line to another line must belong to a given range, as imposed by Equation (26).

$$F_{i,l,t} - W_{i,l,t} \leq vl_l(1 - u_{i,l,t}), \quad \forall\, i \in I_{l-1}, l(l > 1), t(t \neq |T|) \quad (24)$$

$$F_{i,l-1,t} \geq u_{i,l,t}vl_{l-1}, \quad \forall\, i \in I_{l-1}, l(l > 1), t(t \neq |T|) \quad (25)$$

$$S_{\min}\, u_{i,l,t} \leq S_{i,l,t} \leq s_{\max}\, u_{i,l,t}, \quad \forall\, i \in I_{l-1}, l\quad (l > 1), t(t \neq |T|) \quad (26)$$

$$\sum_{i \in I_{l-1}} S_{i,l,t} \leq s_{\max}\left(1 - \sum_{i \in I_l^{Input}} w_{i,l,t}\right), \quad \forall\, l\, (l > 1), t(t \neq |T|) \quad (27)$$

$$W_{i,l,t} = W_{i,l,t-1} + R_{i,l,t} + S_{i,l,t} - S_{i,l+1,t} - \sum_{d \in D_l} D_{i,d,l,t}, \quad \forall\, i \in I_l, l, t > 1 \quad (28)$$

$$W_{i,l,t} = w_{i,l}^{Initial}, \quad \forall i \in I_l^{old}, l, t = 1 \quad (29)$$

$$\sum_{i \in I_l} W_{i,l,t} = vl_l, \quad \forall l, t\, (t \neq |T|) \quad (30)$$

$$\sum_{i \in I_{l-1}} S_{i,l,t} + \sum_{i \in I_l^{Input}} R_{i,l,t} = \sum_{d \in D_l} \sum_{i \in I_{d,l}} D_{i,d,l,t} + \sum_{i \in I_{l-1}} S_{i,l+1,t}, \quad \forall\, l, t(t \neq |T|) \quad (31)$$

$$\min z = \sum_{\substack{t \\ t \neq |T|}} \sum_l \sum_p \sum_{i \in I_l^{Input}} cp_{p,l} RP_{i,p,l,t} + \sum_l \sum_d \sum_p cb_{p,d,l} Border_{p,d,l} + \sum_i IC_i \quad (32)$$

$$IC_i \geq Interface_{p,p'}(y_{i,p} + y_{i-1,p'} - 1), \quad \forall i > 1 \quad (33)$$

**4.6. Tracking Batch Size Over the Time** during time slot $t$, the size of a batch inside line $l$ changes when there are flows going in (due to injection by input node $l$ and material coming from line $l - 1$) or out (due to depots and line $l + 1$). The mass balance in Equation (28) computes the volume of batch $i$ in line $l$ at time point $t$. The initial size of old batch $i$ in every line $l$ is a known parameter, $w_{i,l}^{Initial}$ (Equation (29)).

**4.7. Pipeline Volumetric Balance** at any time, the pipeline is full of product, and thus the sum of batch volumes all together inside each line should be equal to the volume of that line, as imposed by Equation (30).

Moreover, the amount of all batches entered a line $l$ during a time slot should be equal to the amount that leaves this line, as stated in Equation (31).

**4.8. Objective Function** Since pipeline scheduling is a part of oil products supply chain, the aim should be fulfilling demands while minimizing the operational costs as indicated by Equation (32), where the continuous variable $IC_i$ accounts for the interface cost between batches and satisfies Equation (33).

**5. RESULTS AND DISCUSSIONS**

**5.1. Example 1.** This example is taken from [6], considering a pipeline network of two input nodes L1-L2 and three output nodes D1-D3, as shown in the first line of Figure 4. The pipeline transports three products A, B and C at a flow rate belonging to the interval [0.8,1.2] m$^3$/unit.

**TABLE 2.** Product supplies and demands for Example 1

|   | Supplies | | Demands | | |
|---|---|---|---|---|---|
|   | L1 | L2 | D1 | D2 | D3 |
| A | 50 | 20 | 60 | 60 | 0 |
| B | 80 | 60 | 0 | 0 | 100 |
| C | 30 | 40 | 0 | 60 | 0 |

**TABLE 3.** Interface costs (in 100 $)

| Predecessors | Successors | | |
|---|---|---|---|
|   | A | B | C |
| A | - | 22 | 35 |
| B | 24 | - | 21 |
| C | 30 | 32 | - |

**TABLE 4.** Product dependent pumping costs (100 $ per unit)

| Product | Refinery | |
|---|---|---|
|   | Source 1 | Source 2 |
| A | 29.00 | 14.50 |
| B | 34.00 | 17 |
| C | 49.00 | 24.5 |

**TABLE 5.** Computational results for Example 1

| Model type | Optimal solution | CPU time |
|---|---|---|
| Continuous | 812000 | 275.32 |
| Discrete | 805300 | 2355.781 |

The problem goal is to fulfill product demands during a scheduling horizon of 10 days. From Figure 4, at the start of the time horizon, the pipeline includes a sequence of four batches I1 (B), I2 (A), I4 (B), I5 (A) with volumes 20, 30, 10, and 20 units, respectively. Note that batch I3 between two batches I2 and I4 is an empty batch intended for possible injection by middle refinery. The product inventories at the input nodes L1 and L2, and the product demands at the output nodes are illustrated in Table 2 while Table 3 shows the interface reprocessing costs between two products, while Table 4 contains the pumping costs. To solve this example, we select the length of time slots to be equal to 10 h. As the model in [6], we assume that during each time slot, (1) the size of injections at the input nodes should not be greater than 40 units and less than 10 units, and (2) the minimum delivery volume to an active output node can not be less than 5 units.

The best schedule found using the proposed model is given in Table 5 and Figure 4, where the pumping operations at the input nodes are depicted by downward arrows whereas the delivery operations to the output nodes are determined by upward arrows. According to Figure 4, the input node L1 is the only active input node along the pipeline during time slot 1 (i.e., time interval [0,10 h]), injecting 10 units of product A into batch I5 to deliver the same amount of product A in batch I2 to the output node D2. During time interval [10, 50] h (time slots 2-5), 40 units of batch I5 conveying product A are injected by the input node L1. During the injection of batch I5 from L1 taking place in [10, 50] h, the output node D1 receives 40 units of product A contained in batch I5. On the other hand, the input node L2 inserts 20 units of product A into the batch I2 in time slots T2 and T3, and the output node D2 receives 20 inputted units of batch I2. The input node L2 continues pumping the products during time slots 4 to 7, injecting 40 units of product B into batch I4. During time slots 4, 5, output node D2 receives the total volume of batch I2. After that, D3 takes the entire content of batch I1 during [60,80] h (time slots 6, 7).

In time slots 8, 9, and 10 (i.e., during time interval [80, 110] h), L1 inserts a new batch I7(C) with volume 30 units. During this operation, D3 receives 30 units of product B, from batch I4. Again, L1 continues pumping operation in T11 but this time by injecting 10 units of product B in batch I8. At this pumping operation, D2 receives 10 units of batch I5. During time interval [120, 130] h (time slot 12), L2 injects ten units of product B into batch I6. When I6 is being injected, D3 receives ten units of batch I4. At time point T13, the input node, L2 starts to enlarge the size of I7, resulting in transferring the entire volume of batch I4 to D3. The next pumping run is executed at the origin of the pipeline during time slots

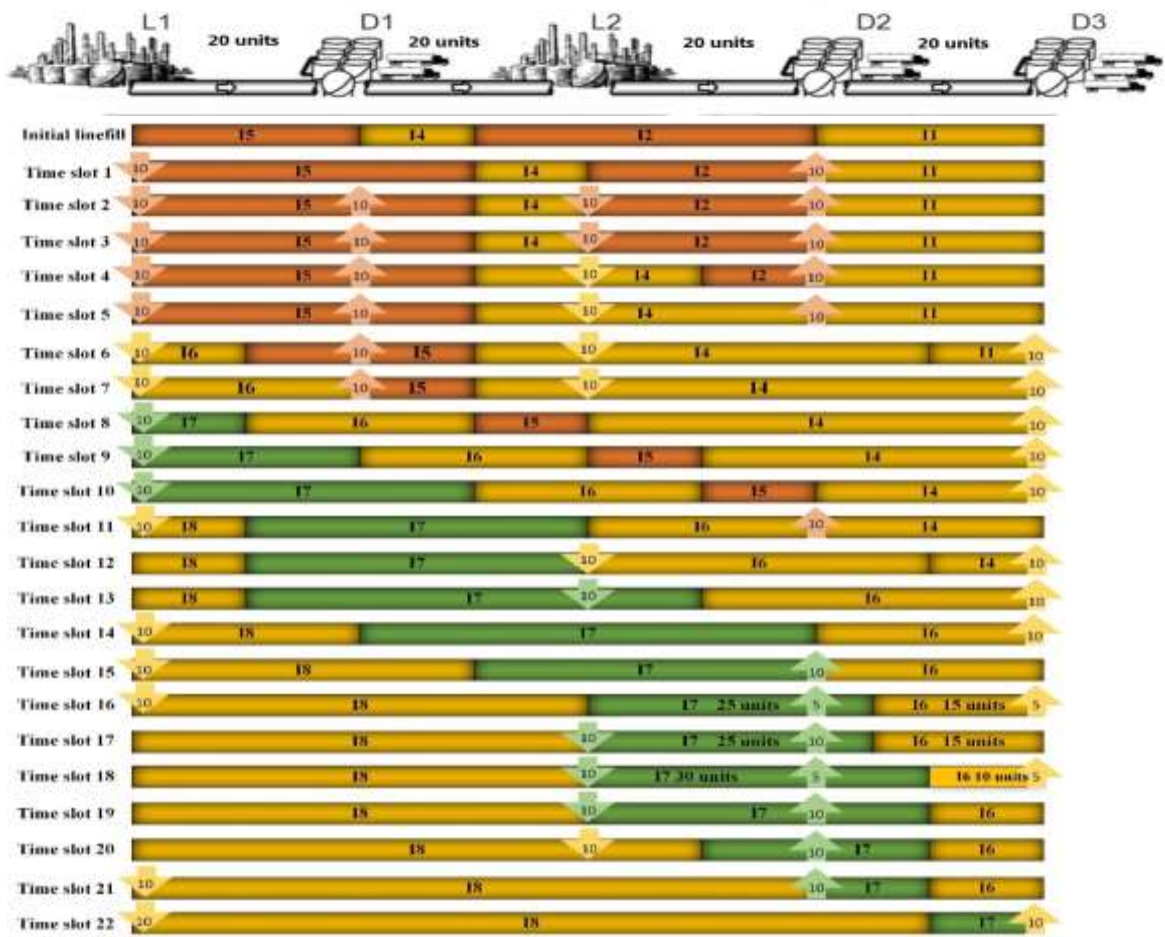

**FIGURE 4.** The best scheduling for Example 1 using the proposed discrete-time

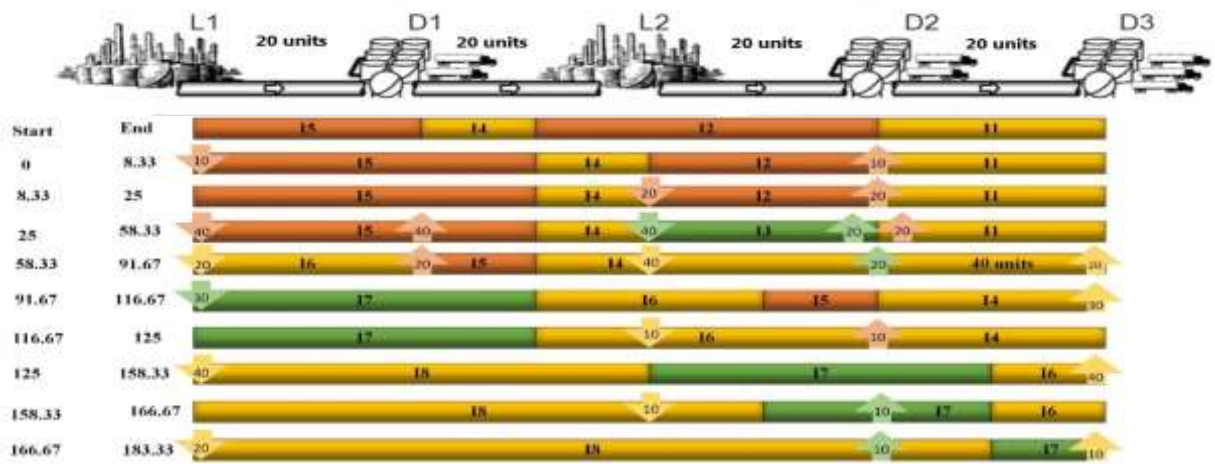

**FIGURE 5.** The best scheduling for Example 1 using the continuous-time method in [6]

T14, T15, and T16, injecting 30 units of batch I8 (B). During this pumping run, the following delivery operations are scheduled: (1) during the time slot T14, the output node D3 receives 10 units of I6, then D2 takes delivery of I7 of the size 10 units in time slot T15 and finally, in time slot 16, (2) D2 withdraws 10 units of I7 during the time slot T15, and (3) at T16, two delivery operations take place in D2 and D3, delivering 5 units of I7 to D2 and 5 units of I6 to D3.

During time interval [170, 200] h, L2 injects product C into the pipeline and then, 30 units of product C to the batch I7. During this operation, 25 units of the batch I& are received by D2 at the time slots T17, T18 and T19. In the middle of the time slot T18, the delivery action of D2 is interrupted and D3 starts to extract five units of I6 during time interval [185, 190] h. In the time slot T20, L2 pumps ten units of the product B inside batch I8 and, again D2 acts as an active delivery node and receives 10 units of batch I7. Finally, the last pumping run is taken place at L1 during time slots T21 and T22. It injects 20 units of product B to deliver 10 units of batch I7 to D2 and 10 units of batch I6 to D3

Figure 5 shows the best schedule for Example 1 using a continuous-time model reported in [6]. It can be seen that the continuous-time model requires fewer time slots to generate the optimal solution, but leads to higher operational cost (see Table 5). This is because the model in [6] fills the empty batch I3 with product C3 and generates **two new interfaces** B-C and C-A during the third pumping operation, i.e., during the time window [25, 53.88]h (see the fourth line of Figure 5). The batch I3 remains empty throughout the planning horizon in the solution found using the proposed discrete-time model, and hence no interface is generated between batch I2 and I4.

**5.2. Example 2** This example is a variant of the first example and considers a time horizon of 5 days. The product demands and the product inventories at the pipeline nodes are given in Table 6. To make a fair comparison between the results provided using the proposed model and those published in [7], we assume that only one of the refineries L1 and L2 can be active during a time window.

To solve this example, we assume that each time slot (window) is composed of 8.5 hours. The computational result for this example is given in Table 7 while Figure 7 depicts the optimal injection and delivery plans during the planning horizon. For example, in Figure 7, during the first two time slots 1 and 2 (i.e., the time interval [0, 17] h), (1) the input node L1 injects 20 units of product A. (2) the output node D2 receives 10 units of product A from batch I2 and the same amount of product A is transferred to the input node D1 from batch I5. By analyzing the results in Table 7 and Figures 6-7, the following conclusions can be summarized: (1) the pipeline operational cost in the proposed model is reduced by 1.5 % compared to that in the continuous-time model in [7], and (2) the optimal solution generated by the proposed model involves 5 product interfaces compared to 7 interfaces in the continuous-time model in [7].

**TABLE 6.** Product supplies and demands for Example 2

| | Supplies | | Demands | | |
|---|---|---|---|---|---|
| | L1 | L2 | D1 | D2 | D3 |
| A | 20 | 10 | 30 | 30 | 0 |
| B | 40 | 30 | 0 | 0 | 50 |
| C | 20 | 20 | 0 | 30 | 0 |

**TABLE 7.** Computational results for Example 2

| Model type | Optimal solution | CPU time |
|---|---|---|
| Continuous | 426700 | 45.3 |
| Discrete | 420200 | 19.156 |

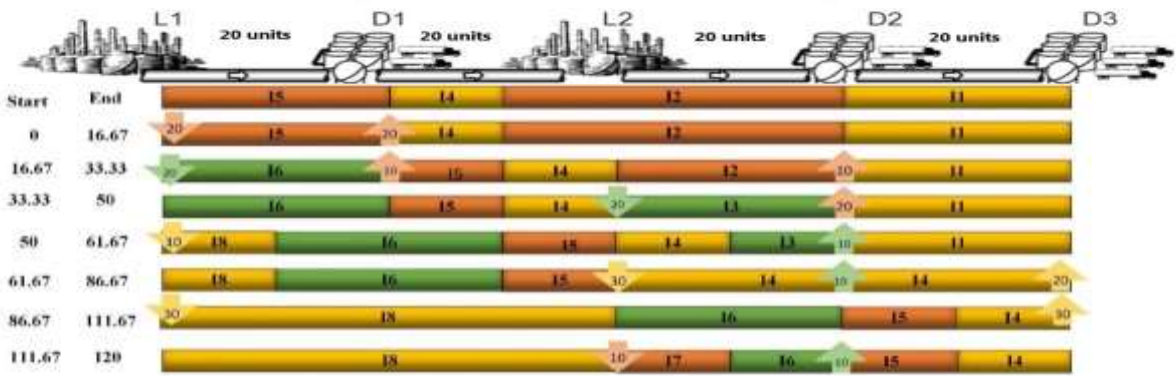

**FIGURE 6.** The best scheduling for Example 2 using the continuous-time method in [7]

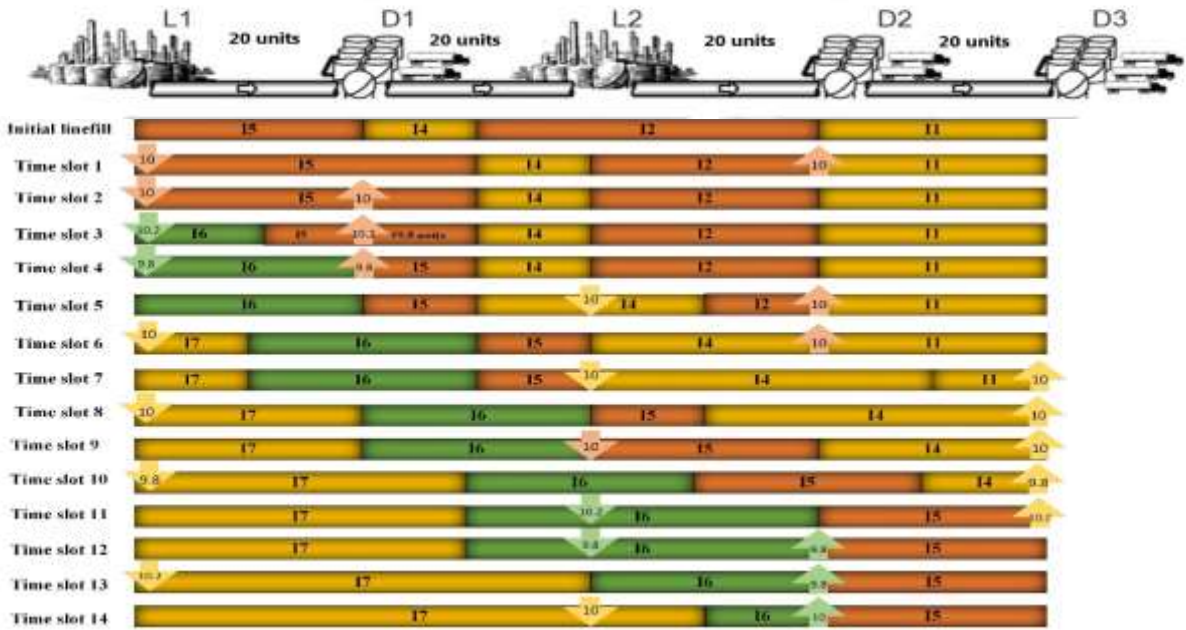

**FIGURE 7.** The best scheduling for Example 2 using the proposed discrete-time

**Table 8.** The effect of backorder cost on optimal solution

| Penalty Cost | Backorder Volume | Objective Value |
|---|---|---|
| 1000 | 126.4 | 172840 |
| 2000 | 83.2 | 287620 |
| 3000 | 60 | 364000 |
| 4000 | 20 | 402200 |
| 4200 | 20 | 406200 |
| 4400 | 20 | 410200 |
| 4600 | 20 | 414200 |
| 5000 | 0 | 420200 |
| 5200 | 0 | 420200 |
| 5400 | 0 | 420200 |

Table 8 and Figure 8 show the effect of backorder penalty cost on the objective function. It can be seen that when the backorder cost coefficient takes a value between [1000, 5000], product demands can not be fully satisfied. For example, when the backorder coefficient is 4200 $/unit, the optimal solution is 402200, and the model fails to satisfy 20 units of the product demands. When

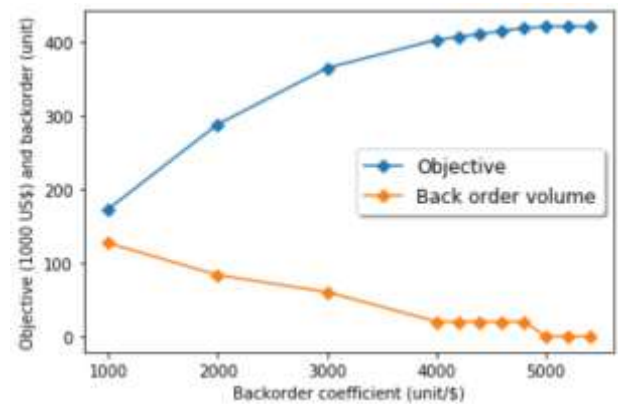

**FIGURE 8.** The effect of backorder cost on the optimal solution

choosing the back order coefficient in [5000,5400], the objective function value does not change and the demands are fully satisfied. Since the main challenge in the pipeline scheduling is to fully meet product demands requested by the output nodes on time, the best value for the backorder coefficient in the objective function should be greater than or equal to 5000 $/unit. The reason why the proposed model generates good solutions when backorder cost belongs [1000, 5000] is that the lower amount of products is pumped into the pipeline (Table 8 comprises of three elements: (1) penalty cost (2)

backorder volume (3) obj. value (objective value). It shows that for example if penalty cost is 4000 $ then the model decides to have 20 units backorder and the objective function will be 402200).

## 6. CONCLUSIONS

This paper addressed scheduling of multiproduct pipelines comprised of multiple refineries and depots, proposing a new discrete-time MILP formulation to model the stated pipeline. The problem goal was to fulfill product demands at the output nodes at minimum costs including pumping, backorder, and interface reprocessing costs. The proposed model can determine decisions related to the pumping and delivery operations at input and output nodes, and to track the size and location of product batches in the pipeline at any time. Moreover, it allows to (1) perform simultaneous pumping operations at input nodes, and (2) to handle the pump and delivery operations at dual-purpose nodes, taking place at the same time window. Two examples were used for testing and validating the proposed approach. The results illustrated that the proposed discrete-time model can lead to better pipeline schedules compared to a similar continuous-time model in the literature. Due to the proper sequencing of product batches in the pipeline, the proposed model was able to notably reduce the interface reprocessing cost, and thus leading to 0.8% and 1.5 % reductions in the operational costs in examples 1 and 2, respectively. The major problem assumption was to assume that all the model parameters as deterministic data. Therefore, future research is needed to consider uncertainty in the model parameters e.g., demand and interface volume. The proposed model can also be extended to consider product demands on multiple intermediate due dates.